\documentclass{article}

\usepackage[latin1]{inputenc}
\usepackage{subfigure}
\usepackage[margin=1in]{geometry}
\usepackage{graphicx}             
\usepackage{amsmath}               
\usepackage{amsfonts}             
\usepackage{amsthm}   
\usepackage{algorithmic}
\usepackage{algorithm}

\newtheorem{thm}{Theorem}[]

\newcommand{\lc}{\left\{}  
\newcommand{\rc}{\right\}}
\newcommand{\ld}{\left[}  
\newcommand{\rd}{\right]}
\newcommand{\lb}{\left(}  
\newcommand{\rb}{\right)} 
\newcommand{\E}{\mathop{\bf E\/}}
\newcommand{\Prob}{\mathop{\bf P\/}}

\newcommand{\remove}[1]{}

\date{}

\title{On the degree distribution of a growing network model}
\author{Linda Farczadi\\ University of Waterloo\\
lindafarczadi@gmail.com \and Nicholas Wormald\thanks{Supported by an ARC Australian Laureate Fellowship. Research supported partly by NSERC}\\ Monash University\\ nick.wormald@monash.edu }
\begin{document}
\maketitle

\begin{abstract}
In this note we make some specific observations on the distribution of the degree of a given vertex in certain model of randomly growing networks. The rule for network growth is the following. Starting with an initial graph of minimum degree at least $k$, new vertices are added one by one. Each new vertex $v$ first chooses a random vertex $w$ to join to, where the probability of choosing $w$ is proportional to its degree. Then $k$ edges are added from $v$ to randomly chosen neighbours of $w$.
\end{abstract}


\section{Introduction}

In this note we make some specific observations on the distribution of the degree of a given vertex in certain model of randomly growing networks. Fix an integer $k \geq 2$. We start with a seed graph $G^1_k$ consisting of one vertex $v_1$ with $k$ loops. For $t \geq 2$ given $G^{t-1}_k$ we obtain $G^{t}_k$ as follows:
\begin{itemize}
\item we add a new vertex $v_{t}$ which connects first to an existing vertex $v_i$ chosen by preferential attachment and then to $k-1$ neighbours of $v_i$ chosen uniformly at random.
\end{itemize}
\noindent We let $V_t$ and $E_t$ denote the vertex set and edge set of $G^t_k$.  Note that  $|V_t| = t$  and $|E_t| = kt$. We denote by $d_{t}(v_i)$ the degree of vertex $v_i$ in $G^t_k$.  We define $ N_{t}(v_i)$ to be the set of neighbours of vertex $v_i$ in $G^t_k$.

We begin with a simple derivation of the expected value of $d_{t}(v_i)$ in Section~\ref{sec: 1.1}, then describe some closely related existing results, and apply them to get a more precise description of   the distribution of $d_{t}(v_i)$ in Section~\ref{sec: 1.4}. 

It is straightforward to modify our results for any given initial seed graph.

\section{Degree distribution}

\subsection{Expected degree of a given vertex}\label{sec: 1.1}

\noindent Fix vertex $v_i$. We want to study $d_n(v_i)$, the degree of vertex $i$ at the $n^{th}$ step in the process. For $t \geq i$ we have 
\begin{align}
\Prob \lb v_{t+1} \text{ connects to $v_i$} | d_t(v_i) \rb &=  \dfrac{ d_{t}(v_i)  }{   2 |E_t| } + \sum_{v_j \in N_{t}(v_i) } \dfrac{ d_{t}(v_j)  }{  2 |E_t| }  \dfrac{ k }{  d_{t}(v_j) } \notag \\
&= \dfrac{ d_{t}(v_i)  }{   2 kt } + d_t(v_i) \lb  \dfrac{ k-1  }{  2kt } \rb \notag \\
&=  \dfrac{ d_{t}(v_i)  } { 2t }. \notag
\end{align}
\noindent Taking expectations of both sides gives
\begin{align}
\Prob \lb v_{t+1} \text{ connects to $v_i$}  \rb &=    \dfrac{  \E \ld d_{t}(v_i) \rd  } { 2t }. \notag
\end{align}
\noindent We then have
\begin{align}
\E \ld d_{t+1}(v_i) \rd  &=    \E \ld d_{t}(v_i) \rd + \Prob \lb v_{t+1} \text{ connects to $v_i$}  \rb \notag \\
&= \E \ld d_{t}(v_i) \rd + \dfrac{  \E \ld d_{t}(v_i) \rd  } { 2t } \notag \\
&= \lb 1 + \frac{1}{2t} \rb \E \ld d_{t}(v_i) \rd. \notag
\end{align}
\noindent Since each vertex has degree $k$ when it joins the graph we have $\E \ld d_{i}(v_i) \rd = k$. We obtain for $1 \leq i \leq n$
\begin{align}
\E \ld d_{n}(v_i) \rd  &=  k \prod_{t=i}^{n-1} \lb 1 + \frac{1}{2t} \rb  \notag \\
&= \dfrac{k \Gamma(i) \Gamma(n+1/2) }{\Gamma(n) \Gamma(i+n/2) } \notag \\
&= k \sqrt{n/i} \lb 1 + O(1/i) \rb . \notag
\end{align}

\subsection{LCD model of Bollob{\'a}s and Riordan}\label{sec: 1.2}

\noindent The LCD model of Bollob{\'a}s and Riordan can be described as follow: start with $G_1^1$ the graph with one vertex and one loop; for $t \geq 2$ given $G_1^{t-1}$ obtain $G_1^t$ by adding one vertex $v_t$ and one edge connecting $v_t$ to an existing vertex $v_i$ chosen randomly with probability given by
\begin{align}
\Prob \lb v_i=s \rb  \begin{cases} \frac{d_{t-1}(s)}{2t-1} & \text{ if $1 \leq s \leq t-1$} \\ \frac{1}{2t-1} & \text{ if $s=t$ } \end{cases} \notag
\end{align}
\noindent Then for a given parameter $k >1$ obtain $G_k^n$  by first constructing $G_1^{kn}$ on vertices $v'_1, v'_2, \cdots, v'_{kn}$ using the process described above. Then  identify vertices $v'_1, \cdots, v'_k$ to form vertex $v_1$ of $G_k^n$, vertices $v'_{k+1}, \cdots, v'_{2k}$ to form vertex $v_2$, and so on. 

\noindent We observe that both the $k$-neighbour model and the Bollob{\'a}s-Riordan model satisfy the following condition 
\begin{align}
\Prob \lb v_{t+1} \text{ connects to $v_i$} | d_t(v_i) \rb &=    \dfrac{ d_{t}(v_i)  } { \sum_{j=1}^{t} d_t(v_j)}. \notag
\end{align}
\noindent This is knows as the Barab{\'a}si-Albert (BA) description.  Hence the degree of a given vertex has the same distribution in both models. In particular the following result concerning the degree sequence of Bollob{\'a}s and Riordan \cite{bollobas2001}, \cite{bollobas2003} applies to the $k$-neighbour model as well.

\begin{thm}  Let $N_n(d)$ be the number of vertices of degree $d$ in $G^n_k$ and define
\begin{align}
\alpha({k,d}) =  \frac{ 2k(k+1)}{d(d+1)(d+2)}. \notag
\end{align}
Then for a fixed $\epsilon > 0$ and  $0 \leq d \leq n^{1/15}$ the following holds with high probability
\begin{align}
(1-\epsilon) \alpha({k,d}) \leq N_n(d) \leq (1+\epsilon) \alpha({k,d}). \notag
\end{align}
\end{thm}

\subsection{General preferential attachment models of Ostroumova et al.}\label{sec: 1.3}

\noindent We can obtain results about the $k$ neighbour model by observing that it belongs to a certain class of general preferential attachment models. Specifically Ostroumova et al.  \cite{ostroumova2012} define the PA-class by considering all random graph models $\mathcal{G}^n_k$ that fit the following description:
\begin{itemize}
\item $G^n_k$ is a graph with $n$ vertices and $kn$ edges obtained from the following random graph process: start at time $n_0$ with an arbitrary seed graph $G^0_k$ with $n_0$ vertices and $kn_0$ edges; at time $t$ obtain the graph $G^{t}_k$ from $G^{t-1}_k$ by adding a new vertex and $k$ edges connecting this vertex to some $k$ vertices of $G^{t-1}_k$.
\end{itemize}
Then $\mathcal{G}^n_k$ belongs to the class PA-class if it satisfies the following conditions for some constants $A$ and $B$:
\begin{align}
\Prob \lb d_{t+1}(v_i) = d_{t}(v_i) | G^t_k \rb &=  1- A \dfrac{d_t(v_i)}{n} - B\dfrac{1}{n} + O \lb \dfrac{ (d_t(v_i))^2}{n^2} \rb 
\end{align}
\begin{align}
\Prob \lb d_{t+1}(v_i) = d_{t}(v_i) +1 | G^t_k \rb &=  A \dfrac{d_t(v_i)}{n} + B\dfrac{1}{n} + O \lb \dfrac{ (d_t(v_i))^2}{n^2} \rb 
\end{align}
\begin{align}
\Prob \lb d_{t+1}(v_i) = d_{t}(v_i) +j | G^t_k \rb &=  O \lb \dfrac{ (d_t(v_i))^2}{n^2} \rb \indent 2 \leq j \leq k
\end{align}
\begin{align}
\Prob \lb d_{t+1}(v_{t+1}) = k +j | G^t_k \rb &=  O \lb \dfrac{ 1}{n} \rb \indent 1 \leq j \leq k
\end{align}
\noindent Then we can observe that our model belongs to this PA-class with parameters $A = 1/2$ and $B=0$.  Then the following two results from \cite{ostroumova2012} apply to our model.

\begin{thm} Let $N_n(d)$ be the number of vertices of degree $d$ in $G^n_k$ and $\theta(X)$ be an arbitrary function such that $| \theta(X) | < X$. There exists a constant $C > 0$ such that for any $d \geq k$ we have
\begin{align}
\E \ld N_n(d) \rd = \alpha(k,d) \lb n + \theta(C d^4) \rb \notag
\end{align}
where
\begin{align}
\alpha(k,d)= \frac{ 2 k(k+1) }{ d(d+1)(d+2) } \sim 2 k(k+1) d^{-3}.  \notag
\end{align}
\end{thm}

\begin{thm} For any $\delta > 0$ there exists a function $\psi(n) = o(n)$ such that for any $k \leq  d \leq n^{\frac{1}{8} - \frac{\delta}{4}}$
\begin{align}
\lim_{n \rightarrow \infty} \Prob \lb \left| N_n(d) - \E \ld N_n(d) \rd \right| \geq \frac{ \psi(n)}{d^3} \rb = 0. \notag
\end{align}
\end{thm}

\subsection{Urn models}\label{sec: 1.4}

\noindent We can obtain the distribution of $d_n(v_i)$ by using an urn model. Our urn contains balls of two colours: white and black. White balls represent edge-ends incident with vertex $i$ and black balls represent edge-ends not incident with vertex $i$. Suppose that the urn initially has $a_0$ white balls and $b_0$ black balls where
\begin{align}
a_0 &= k \notag \\
b_0 &= 2i-k \notag \\
t_0 &= a_0 + b_0 = 2i \notag.
\end{align}
At each step, one ball is drawn randomly from the urn.  If the drawn ball is white, replace it and put an additional   $\alpha$ white and $\sigma-\alpha$ black. If it is black, replace and put $\sigma$ more black balls. 

We now introduce some relevant results about urn models from Flajolet et al. \cite{flajolet}. Consider  a triangular urn with replacement matrix
\begin{align}
\begin{pmatrix}
\alpha & \sigma-\alpha \\
0 & \sigma
\end{pmatrix} \notag
\end{align}
Let  $H_n \binom{a_0 \, a}{b_0 \, b}$ be the number of histories of length $n$ that start in configuration $(a_0,b_0)$ and end in configuration $(a,b)$. Note that $a+b = a_0 + b_0 + \sigma n$. Then the generating function of urn histories is defined as
\begin{align}
H(x,1,z) := \sum_{n,a} H_n \binom{a_0 \, a}{b_0 \, b} x^a \dfrac{z^n}{n!}  \notag
\end{align}
and is given by
\begin{align}
H(x,1,z) &= x^{a_0} \lb 1- \sigma z \rb^{-b_0/\sigma} \lb 1 - x^\alpha \lb 1- \lb  1- \sigma z \rb^{\alpha/\sigma} \rb \rb^{-a_0/\alpha}. \notag 
\end{align}
\noindent Letting $\Delta := \lb  1- \sigma z \rb^{-1/\sigma}$ we have from~\cite[Equation (74)]{flajolet}
\begin{align}
H(x,1,z) &= x^{a_0} \Delta^{b_0} \lb 1 - x^\alpha \lb 1- \Delta^{-\alpha} \rb \rb^{-a_0/\alpha} \notag 
\end{align}
(The reader may notice that~\cite{flajolet} is not entirely consistent, but in that paper, the `balls of the first type' do always correspond to the first row of the replacement matrix.)  Now let $A_m$ be the number of white balls in the urn after $m$ trials.  The probability that this equals $a_0 + x \alpha$   for some $0 \leq x \leq m$  is given by~\cite[Equation (75)]{flajolet} as
\begin{align}
\Prob \lb A_m = a_0 + x \alpha \rb &=    \binom{x+ \frac{a_0}{\alpha}-1}{x} \sum_{i=0}^{x} (-1)^i \binom{x}{i} 
\frac{ [m+(b_0-\alpha i)/\sigma -1 ]_m}
{[m+ t_0/\sigma  -1]_m}  
\notag
\end{align}
where $[\cdot]_m$ denotes falling factorial.
 
We next explain why the urn results apply to the random network process. At any given step, let $s$ denote the total degrees of the vertices. If a vertex $v$ has degree $d$ then the probability it is chosen as the first vertex is $d/s$. The probability it is chosen as the second vertex via any given one of its $d$ incident edges is $(k-1)/s$. Hence, the probability it receives a new edge is $dk/s$. (This is similar to the derivation in Section 2.1.) So we may use an urn with $kd$ white balls and $s$ black balls. At each step, $s$ increases by $2k$, whilst $kd$ increases by $k$ if $v$ receives a new edge, and by 0 otherwise. Hence, at each step the number of white balls is $k$ times the degree of $v$, and the parameters are $\alpha=k$ and $\sigma=2k$. The initial number of white balls, $a_0$, is $k$ times the initial degree of the vertex, and the initial number of black balls is $2k$ times the initial number of edges.

Suppose the initial graph is a copy of $K_j$. Then for any of the $j$ initial vertices we have $\alpha=k$, $a_0=k(j-1)$, $t_0=j(j-1)$ and $b_0=t_0-a_0= (j-k)(j-1)$. Also $\sigma=2k$. For each of those $j$ vertices $v_i$ (which initially have degree $j-1$), there are $m=n-j$ trials. So the probability that the corresponding urn process finishes with $k(j-1 + x )$ balls is
$$
\Prob \lb d_n(v_i) = j-1 + x \rb  =  \binom{x+j-2}{x} \sum_{u=0}^{x} (-1)^u \binom{x}{u}
\frac{\left[n-j+\frac{(j-1)^2-ku}{2k}-1\right]_{n-j}  }
{[n-j+R -1]_{n-j}}
$$
where $R=j(j-1)/(2k)$. On the other hand, if $i>j$ then $t_0=j(j-1)+2(i-j)k$, $a_0=k^2$, $b_0=t_0-a_0$, $m=n-i$ and so
$$
\Prob \lb d_n(v_i) = k + x \rb  = 
 \binom{x+k-1}{x} \sum_{u=0}^{x} (-1)^u \binom{x}{u}\frac{
 \left[n -j+R- (u+k)/2 -1\right]_{n-i}}
 {  [n-j +R-1]_{n-i}}.  
$$

\noindent We can also obtain the moments of $A_n$ directly from the generating function of urn histories. In particular, as explained in~\cite{flajolet}, the first and second moments are given by
\begin{align}
\E \ld A_n \rd &= \frac{ \Gamma(n+1) \Gamma \lb \frac{t_0}{\sigma} \rb }{ \sigma^n \Gamma  \lb \frac{t_0}{\sigma} +n \rb} \ld z^n \rd \lb \frac{ \partial H(x,1,z)}{ \partial x} \rb_{x=1} \notag \\
&=  \frac{ \Gamma(n+1) \Gamma \lb \frac{t_0}{\sigma} \rb }{ \sigma^n \Gamma  \lb \frac{t_0}{\sigma} +n \rb} \ld z^n \rd a_0 \Delta^{t_0 + \alpha}  \notag \\
\E \ld A_n^2 \rd &= \frac{ \Gamma(n+1) \Gamma \lb \frac{t_0}{\sigma} \rb }{ \sigma^n \Gamma  \lb \frac{t_0}{\sigma} +n \rb} \ld z^n \rd \lb \frac{ \partial^2 H(x,1,z)}{ \partial x^2} \rb_{x=1} \notag \\
&=  \frac{ \Gamma(n+1) \Gamma \lb \frac{t_0}{\sigma} \rb }{ \sigma^n \Gamma  \lb \frac{t_0}{\sigma} +n \rb} \ld z^n \rd \lb a_0 \lb a_0+\alpha \rb \Delta^{t_0+2\alpha}- a_0(\alpha+1) \Delta^{t_0+\alpha} \rb. \notag
\end{align}
\noindent Performing coefficient extraction gives
\begin{align}
\ld z^n \rd a_0 \Delta^{t_0 + \alpha} &=  a_0 \ld z^n \rd \lb 1- \sigma z \rb^{-\frac{t_0 + \alpha}{\sigma}} \notag \\
&=a_0 \binom{ -\frac{t_0 + \alpha}{\sigma}}{n} \lb -\sigma \rb ^n \notag \quad \text{ using $\ld t ^n\rd (1+at)^r = \binom{r}{n} a^n$ } \notag \\
&= a_0 \binom{ n+\frac{t_0 + \alpha}{\sigma}-1}{n}  \sigma^n \quad \text{ using $\binom{-r}{n} = \binom{n+r-1}{n} (-1)^n$ } \notag \\
&= a_0 \sigma^n \frac{ \Gamma \lb  n+\frac{t_0 + \alpha}{\sigma} \rb }{ \Gamma \lb n +1 \rb \Gamma \lb \frac{t_0 + \alpha}{\sigma} \rb} \quad \text{ using $\binom{x}{y} = \frac{ \Gamma \lb  x+1  \rb }{ \Gamma \lb y +1 \rb \Gamma \lb x-y +1 \rb}$ } \notag 
\end{align}
\noindent And by linearity:
\begin{align}
\ld z^n \rd \lb a_0 \lb a_0+\alpha \rb \Delta^{t_0+2\alpha}- a_0(\alpha+1) \Delta^{t_0+\alpha} \rb &= a_0 \lb a_0+\alpha \rb \ld z^n \rd \Delta^{t_0+2\alpha} - a_0(\alpha+1) \ld z^n \rd \Delta^{t_0+\alpha} \notag \\
&=  a_0  \ld \lb a_0+\alpha \rb \sigma^n \frac{ \Gamma \lb  n+\frac{t_0 + 2\alpha}{\sigma} \rb }{ \Gamma \lb n +1 \rb \Gamma \lb \frac{t_0 + 2\alpha}{\sigma} \rb} - (\alpha+1) \sigma^n \frac{ \Gamma \lb  n+\frac{t_0 + \alpha}{\sigma} \rb }{ \Gamma \lb n +1 \rb \Gamma \lb \frac{t_0 + \alpha}{\sigma} \rb} \rd. \notag 
\end{align}
\noindent Plugging these coefficients back in the equations for the first and second moment we obtain
\begin{align}
\E \ld A_n \rd &=  \frac{ a_0 \Gamma \lb \frac{t_0}{\sigma} \rb \Gamma \lb  n+\frac{t_0 + \alpha}{\sigma} \rb }{ \Gamma \lb \frac{t_0 + \alpha}{\sigma} \rb \Gamma  \lb \frac{t_0}{\sigma} +n \rb}  \notag \\
\E \ld A_n^2 \rd &= a_0 \ld \frac{ (a_0 + \alpha) \Gamma \lb \frac{t_0}{\sigma} \rb \Gamma \lb  n+\frac{t_0 + 2\alpha}{\sigma} \rb }{ \Gamma \lb \frac{t_0 + 2\alpha}{\sigma} \rb \Gamma  \lb \frac{t_0}{\sigma} +n \rb} + \frac{ (\alpha+1) \Gamma \lb \frac{t_0}{\sigma} \rb \Gamma \lb  n+\frac{t_0 + \alpha}{\sigma} \rb }{ \Gamma \lb \frac{t_0 + \alpha}{\sigma} \rb \Gamma  \lb \frac{t_0}{\sigma} +n \rb} \rd \notag
\end{align}
\text{}\\
\noindent Applying these results to our example where $a_0 =k$, $b_0 = 2i -k$, $\alpha =1$, $\sigma =2$ and the number of trials is $n-i$ we obtain
\begin{align}
\E \ld d_n(i) \rd &=  \frac{ k \Gamma \lb i \rb \Gamma \lb  n+1/2 \rb }{ \Gamma \lb i+1/2 \rb \Gamma  \lb n \rb}  \notag \\
\E \ld d_{n}(i)^2 \rd &= k \ld \frac{ (k + 1) \Gamma \lb i \rb \Gamma \lb  n+1 \rb }{  \Gamma \lb i+1 \rb \Gamma  \lb n \rb} + \frac{ 2 \Gamma \lb i \rb \Gamma \lb  n+1/2 \rb }{ \Gamma \lb i+1/2 \rb \Gamma  \lb n \rb} \rd \notag \\
&=  \frac{k (k + 1)n }{ i} + 2 \E \ld d_n(i) \rd . \notag
\end{align}
\noindent We note that the value for the first moment matches the one obtained previously by solving the simple recursion in Section \ref{sec: 1.1}.


 \newcommand{\etalchar}[1]{$^{#1}$}

\end{document}